# Small deviations of stable processes and entropy of the associated random operators

FRANK AURZADA[1], MIKHAIL LIFSHITS[2] and WERNER LINDE[3]

[1]*Technische Universität Berlin, Institut für Mathematik, Sekr. MA7-5, Str. des 17, Juni 136, 10623 Berlin, Germany. E-mail: aurzada@math.tu-berlin.de*

[2]*St. Petersburg State University, 198504 Stary Peterhof, Department of Mathematics and Mechanics, Bibliotechnaya pl. 2, Russia. E-mail: lifts@mail.rcom.ru*

[3]*Friedrich-Schiller-Universität Jena, Institut für Stochastik, Ernst-Abbe-Platz 2, 07743 Jena, Germany. E-mail: lindew@minet.uni-jena.de*

We investigate the relation between the small deviation problem for a symmetric $\alpha$-stable random vector in a Banach space and the metric entropy properties of the operator generating it. This generalizes former results due to Li and Linde and to Aurzada. It is shown that this problem is related to the study of the entropy numbers of a certain random operator. In some cases, an interesting gap appears between the entropy of the original operator and that of the random operator generated by it. This phenomenon is studied thoroughly for diagonal operators. Basic ingredients here are techniques related to random partitions of the integers. The main result concerning metric entropy and small deviations allows us to determine or provide new estimates for the small deviation rate for several symmetric $\alpha$-stable random processes, including unbounded Riemann–Liouville processes, weighted Riemann–Liouville processes and the ($d$-dimensional) $\alpha$-stable sheet.

*Keywords:* Gaussian processes; metric entropy; random operators; Riemann–Liouville processes; small deviations; stable processes

## 1. Introduction

Let $[E, \|\cdot\|_E]$ be a (real) Banach space with (topological) dual space $E'$. We endow $E'$ with the weak-$*$-topology and denote by $\mathcal{B}_\sigma(E')$ the corresponding $\sigma$-field. Now, an $E'$-valued random vector on $(\Omega, \mathbb{P})$ is always understood to be measurable with respect to this $\sigma$-field. Such a vector $X$ is said to be *symmetric $\alpha$-stable* (as usual, we write S$\alpha$S for short) for some $\alpha \in (0, 2]$ if there exist a measure space $(S, \sigma)$ and a (linear, bounded) operator $u : E \to L_\alpha(S, \sigma)$ such that

$$\mathbf{E} e^{\mathrm{i}\langle z, X \rangle} = e^{-\|u(z)\|_\alpha^\alpha}, \qquad z \in E. \tag{1.1}$$

In this case, we say that $X$ is generated by the operator $u$. This approach is very useful for investigating symmetric $\alpha$-stable random processes with paths in $E'$. We refer to Section 5 of Li and Linde (2004) or Section 7.1.1 below for a discussion of how all natural examples of S$\alpha$S processes fit into this framework. For example, if $u$ from $L_p[0,1]$ to $L_\alpha[0,1]$ is defined by $(uf)(t) := \int_t^1 f(s)\,ds$, $f \in L_p[0,1]$, then the random vector $Z_\alpha$ generated via (1.1) is nothing but the symmetric $\alpha$-stable Lévy motion (see Samorodnitsky and Taqqu (1994) for the definition) with paths regarded in $L_{p'}[0,1]$ (as usual, we use $p'$ to denote the conjugate of $p$, that is, $1/p + 1/p' = 1$). In particular, for $\alpha = 2$, we consider the Wiener process in $L_{p'}[0,1]$.

A symmetric 2-stable vector is centered Gaussian. In this case, there exist tight relations between the degree of compactness of $u: E \to L_2$ and small deviation properties of the generated random vector $X$.

To make this more precise, let us introduce the small deviation function

$$\phi(X, \varepsilon) := -\log \mathbb{P}(\|X\|_{E'} < \varepsilon) \tag{1.2}$$

of an $E'$-valued random vector $X$. To measure the degree of compactness of the corresponding operator $u$, we use the dyadic entropy numbers defined as follows: if $u$ is a bounded linear operator between the Banach spaces (or more general quasi-normed spaces) $E$ and $F$, then we let

$$e_n(u) := \inf\{\varepsilon > 0 \mid \exists y_1, \ldots, y_{2^{n-1}} \in F, \forall z \in E, \|z\| \leq 1, \exists i \leq 2^{n-1}, \|u(z) - y_i\| \leq \varepsilon\}.$$

As can be easily seen, an operator $u$ is compact if and only if the corresponding entropy numbers tend to zero. Thus, their behavior as $n \to \infty$ describes the degree of compactness of $u$.

Before we state the results, let us establish some more notation. We write $f \preceq g$ or $g \succeq f$ if $\limsup f/g < \infty$, while the equivalence $f \approx g$ means that we have both $f \preceq g$ and $g \preceq f$. Moreover, $f \lesssim g$ or $g \gtrsim f$ indicate that $\limsup f/g \leq 1$. Finally, the strong equivalence $f \sim g$ means that $\lim f/g = 1$.

Using this notation, we can now state the aforementioned relation between properties of $X$ and the generating operator $u$ in the Gaussian case.

**Proposition 1.1 (Kuelbs and Li (1993), Li and Linde (1999)).** *Assume that $X$ is an $E'$-valued Gaussian vector that is generated by the operator $u: E \to \ell_2$. Let $\tau > 0$ and let $L$ be a slowly varying function at infinity such that $L(t) \approx L(t^p)$ for all $p > 0$. The following implications then hold:*

(a) *we have*

$$e_n(u) \succeq n^{-1/2-1/\tau} L(n) \quad \Leftrightarrow \quad \phi(X, \varepsilon) \succeq \varepsilon^{-\tau} L(1/\varepsilon)^\tau,$$

*where, for "$\Leftarrow$", the additional assumption $\phi(X, \varepsilon) \approx \phi(X, 2\varepsilon)$ is required;*

(b) *we have*

$$e_n(u) \preceq n^{-1/2-1/\tau} L(n) \quad \Leftrightarrow \quad \phi(X, \varepsilon) \preceq \varepsilon^{-\tau} L(1/\varepsilon)^\tau.$$

It is natural to ask whether or not these implications can be transferred to the non-Gaussian setup of symmetric $\alpha$-stable vectors. In this case, the following is known.

**Proposition 1.2 (Li and Linde (2004), Aurzada (2007b)).** *Let $X$ be an $E'$-valued symmetric $\alpha$-stable vector generated by an operator $u: E \to L_\alpha(S,\sigma)$. Let $\tau > 0$ and $\theta \in \mathbb{R}$ be given, where, additionally, $\tau < \alpha/(1-\alpha)$ for $0 < \alpha < 1$. Then,*

(a) $e_n(u) \succeq n^{1/\alpha - 1/\tau - 1}(\log n)^{\theta/\tau}$ *implies* $\phi(X,\varepsilon) \succeq \varepsilon^{-\tau}(-\log \varepsilon)^\theta$,
(b) $\phi(X,\varepsilon) \preceq \varepsilon^{-\tau}(-\log \varepsilon)^\theta$ *implies* $e_n(u) \preceq n^{1/\alpha - 1/\tau - 1}(\log n)^{\theta/\tau}$ *and*
(c) *the respective converse in the above implications does not hold in general.*

This result shows that, unfortunately, only two of the four implications from Proposition 1.1 can be transferred to the non-Gaussian case. In particular, probably the most interesting and useful implication (upper estimates for $e_n(u)$ yield those for $\phi(X,\varepsilon)$) is not valid in general. The basic goal of this article is to investigate this implication more thoroughly. It turns out that if we take the entropy numbers of $u: E \to L_\alpha(S,\sigma)$ regarded as operator into $L_\infty(S,\sigma)$, then this implication is also valid. Let us mention (see Section 2 below) that we may always assume that an operator generating an S$\alpha$S vector can be factorized over $L_\infty(S,\sigma)$. Therefore, in all cases of interest, the operator $u_\infty$, which is simply $u$ acting from $E$ to $L_\infty(S,\sigma)$, is well defined. We fix this notation for $u$ and $u_\infty$ throughout the article.

The main result of this paper is the following.

**Theorem 1.3.** *Let a symmetric $\alpha$-stable $E'$-valued vector $X$ be generated by an operator $u: E \to L_\alpha(S,\sigma)$, where $\sigma(S) < \infty$. Suppose that $u$ maps $E$ even into $L_\infty(S,\sigma)$ and that*

$$e_n(u_\infty) = e_n(u: E \to L_\infty) \preceq n^{1/\alpha - 1/\tau - 1} L(n)$$

*for some $\tau > 0$ and some slowly varying function $L$ such that $L(t) \approx L(t^p)$ for all $p > 0$. Then,*

$$\phi(X,\varepsilon) \preceq \varepsilon^{-\tau} L(1/\varepsilon)^\tau.$$

The proof is postponed to Section 4.2.

Note that $e_n(u) \leq \sigma(S)^{1/\alpha} e_n(u_\infty)$. Thus, Theorem 1.3 is weaker than the corresponding result in the Gaussian case. Nevertheless, there are many examples of interest where the entropy numbers of $u$ and $u_\infty$ have the same asymptotic order. Consequently, for those operators, the implication "$\Rightarrow$" in (b) of Proposition 1.2 is also valid. Below, we shall give several examples of this situation.

This article is structured as follows. In Section 2, we analyze decomposed operators from a Banach space $E$ into $L_\alpha(S,\sigma)$, $0 < \alpha < 2$. It is shown that such operators are associated with random operators $v$ mapping $E$ into $\ell_2$. As a consequence, we get the well-known fact that each $E'$-valued symmetric stable vector is a mixture of suitable Gaussian ones. This fact is the basic ingredient of the proof of Theorem 1.3. In Section 3, we investigate compactness properties of the random operator $v$. In particular, we show that the entropy numbers of $u_\infty$ and those of the associated random operator $v$ are closely related. In Section 4, we state and prove a zero–one law for the entropy behavior of the random operator $v$. As a consequence, the entropy numbers of $v$ possess a.s. the same degree of compactness. Furthermore, in that section, Theorem 1.3 is proved.

Although we could shed some light on the relation between the small deviations of S$\alpha$S vectors and the entropy behavior of the generating operator, several interesting questions remain open. The most important ones are presented in Section 5. Besides, an interesting phenomenon is considered: in some cases, a surprising gap appears between the entropy of the original operator and that of its associated random one. In Section 6, this gap is investigated thoroughly in the case of diagonal operators. This problem finally leads to the investigation of diagonal operators with random diagonal. In the authors' opinion, the results of that section could be of independent interest.

Theorem 1.3 gives new bounds or clarifies the small deviation rate for several examples of symmetric stable processes. These examples are considered in Section 7, including unbounded Riemann–Liouville processes, weighted Riemann–Liouville processes and the ($d$-dimensional) $\alpha$-stable sheet. Finally, we give a short and direct proof of a result of M. Ryznar concerning $\alpha$-stable vectors with $0 < \alpha < 1$.

## 2. Representation of decomposed operators

The aim of this section is to analyze the structure of the operator $u$ in (1.1). In particular, it can be decomposed (see below) and is thus associated with a random operator $v$ corresponding to a random Gaussian vector. As a consequence, the stable distribution of the vector $X$ may be represented as a suitable mixture of Gaussian ones.

As before, let $[E, \|\cdot\|_E]$ be some (real) Banach space and let $(S, \sigma)$ be a measure space. An operator $u$ from $E$ into $L_\alpha(S, \sigma)$ for some $\alpha > 0$ is said to be *order bounded* provided there is some function $f \in L_\alpha(S, \sigma)$ such that

$$|(uz)(s)| \leq f(s), \qquad \sigma\text{-a.s. for } z \in E, \|z\| \leq 1.$$

A useful equivalent formulation (see Vakhania *et al.* (1985)) is as follows: there exists a $\mathcal{B}_\sigma(E')$-measurable function $\varphi$ from $S$ into $E'$, the topological dual of $E$, such that

$$\int_S \|\varphi(s)\|_{E'}^\alpha \, \mathrm{d}\sigma(s) < \infty \quad \text{and} \quad u(z) = \langle z, \varphi \rangle, \qquad z \in E. \tag{2.1}$$

Let us say that $\varphi$ decomposes the operator $u$. In particular, whenever $s \in S$ is fixed, for those $u$, the mapping $z \mapsto (uz)(s)$ is a well-defined linear functional on $E$.

We note that the operator $u : E \to L_\alpha(S, \sigma)$ generating an $E'$-valued vector as in (1.1) may always be chosen to be order bounded. This follows from Tortrat's theorem concerning the spectral representation of symmetric stable measures (see Tortrat (1976)). In Section 5.2, we come back to the spectral representation as a natural choice for the generating operator.

Let us also prove that one can always use a bounded decomposing function $\varphi$ and a finite measure space $(S, \sigma)$.

**Proposition 2.1.** *Let $u: E \to L_\alpha(S, \sigma)$ be a decomposed operator. There then exists a finite measure $\tilde{\sigma}$ on $S$ and an operator $\tilde{u}: E \to L_\alpha(S, \tilde{\sigma})$ such that $\tilde{u}$ is decomposed by a function $\tilde{\varphi}$ such that $\|\tilde{\varphi}(s)\| \leq 1$, $s \in S$, and $\|u(z)\|_\alpha = \|\tilde{u}(z)\|_\alpha$ for all $z \in E$.*

**Proof.** Suppose that $\varphi$ decomposes an operator $u$ as in (2.1). Set $\tilde{\varphi}(s) := \varphi(s)/\|\varphi(s)\|_{E'}$ for $s \in S$ and define $\tilde{\sigma}$ on $S$ by $\mathrm{d}\tilde{\sigma}(s) := \|\varphi(s)\|_{E'}^\alpha \, \mathrm{d}\sigma(s)$. By the properties of $\varphi$, this measure is finite. Finally, the operator $\tilde{u}: E \to L_\alpha(S, \tilde{\sigma})$ is given by

$$\tilde{u}(z) := \langle z, \tilde{\varphi} \rangle, \qquad z \in E.$$

Of course, we have $\|u(z)\|_{L_\alpha(S,\sigma)} = \|\tilde{u}(z)\|_{L_\alpha(S,\tilde{\sigma})}$ and this completes the proof. □

**Remark.** Note that $u$ and $\tilde{u}$ possess the same compactness properties. Hence, we may (and will) assume that the decomposing function $\varphi$ of $u$ has the additional property

$$\|\varphi(s)\|_{E'} \leq 1, \qquad s \in S, \tag{2.2}$$

and that the underlying measure $\sigma$ is finite.

The following result from Li and Linde (2004) (Proposition 2.1 there) is crucial for our further investigation. In contrast to Li and Linde (2004), we formulate it directly for operators on $E$ (our $v_\delta$ correspond to $v_\delta^*$ in Li and Linde (2004)).

**Proposition 2.2.** *Suppose that $0 < \alpha < 2$ and let $u: E \to L_\alpha(S, \sigma)$ be order bounded. There then exists a probability space $(\Delta, \mathbb{Q})$ and, for $\delta \in \Delta$, an operator $v_\delta: E \to \ell_2$ such that*

$$\mathrm{e}^{-\|u(z)\|_\alpha^\alpha} = \int_\Delta \mathrm{e}^{-\|v_\delta(z)\|_2^2} \, \mathrm{d}\mathbb{Q}(\delta) = \mathbf{E}_\mathbb{Q} \mathrm{e}^{-\|v(z)\|_2^2} \qquad \textit{for all } z \in E. \tag{2.3}$$

In the last expression of (2.3), we omitted the $\delta$ as it is common for random variables. Here, and in the following, we often write $v$ instead of $v_\delta$. We stress, however, that $v$ denotes a *random* operator. In the same way, we shall also often replace the integral with respect to $\mathbb{Q}$ by $\mathbf{E}_\mathbb{Q}$.

For our further investigation, it is important to have more information about the random operator $v = v_\delta$. For this purpose, choose an i.i.d. sequence $(V_j)_{j \geq 1}$ of $S$-valued random variables with common distribution $\sigma/\sigma(S)$. Furthermore, let $(\zeta_j)_{j \geq 1}$ be an i.i.d. sequence of standard exponential random variables. Define $\Gamma_j$ by

$$\Gamma_j := \zeta_1 + \cdots + \zeta_j$$

and suppose that the $V_j$, as well as the $\zeta_j$, are defined on $(\Delta, \mathbb{Q})$ and that these two sequences are independent. Finally, set

$$c_\alpha := \sqrt{2} \bigg( \int_0^\infty x^{-\alpha} \sin x \, \mathrm{d}x \bigg)^{-1/\alpha} (\mathbf{E}|\xi|^\alpha)^{-1/\alpha},$$

where $\xi$ is standard normal. Then, $v\colon E \to \ell_2$ admits the following representation:

$$v(z) = c_\alpha \sigma(S)^{1/\alpha}((uz)(V_j)\Gamma_j^{-1/\alpha})_{j=1}^\infty, \qquad z \in E. \tag{2.4}$$

Recall that $(uz)(V_j)$ has to be understood as $\langle z, \varphi(V_j) \rangle$, where $\varphi$ is the function decomposing $u$.

Let us now define the (random) operators $w\colon E \to \ell_\infty$ and $D\colon \ell_\infty \to \ell_2$ by

$$w(z) := ((uz)(V_j))_{j=1}^\infty, \qquad z \in E, \tag{2.5}$$

and

$$D(y) := c_\alpha \sigma(S)^{1/\alpha}(\Gamma_j^{-1/\alpha} y_j)_{j=1}^\infty, \qquad y = (y_j)_{j \geq 1} \in \ell_\infty. \tag{2.6}$$

Note that both operators are well defined $\mathbb{Q}$-almost surely. Indeed, if $\varphi$ is the decomposing function of $u$, by (2.2), it follows that

$$|(uz)(V_j)| = |\langle z, \varphi(V_j) \rangle| \leq \|z\|.$$

On the other hand, the strong law of large numbers implies that $\lim_{j \to \infty} \Gamma_j/j = 1$. Thus, if $y = (y_j)_{j \geq 1}$ is in $\ell_\infty$, then the sequence $(\Gamma_j^{-1/\alpha} y_j)_{j \geq 1}$ is almost surely square summable since $0 < \alpha < 2$.

Summarizing the previous remarks, we get the following.

**Proposition 2.3.** *Suppose that $0 < \alpha < 2$ and let $u\colon E \to L_\alpha(S, \sigma)$ be a decomposed operator. We then have*

$$\mathrm{e}^{-\|u(z)\|_\alpha^\alpha} = \mathbf{E}_\mathbb{Q} \mathrm{e}^{-\|v(z)\|_2^2}, \qquad z \in E, \tag{2.7}$$

*where $v = D \circ w$, with $D$ and $w$ defined by (2.6) and (2.5), respectively.*

Let $v\colon E \to \ell_2$ be the operator representing $\|u(z)\|_\alpha$ as in Propositions 2.2 and 2.3. As shown by Sztencel (1984), there exist $E'$-valued (centered Gaussian) random vectors $Y = Y_\delta$, $\delta \in \Delta$, such that, $\mathbb{Q}$-almost surely,

$$\mathbf{E}\mathrm{e}^{\mathrm{i}\langle z, Y \rangle} = \mathrm{e}^{-\|v(z)\|_2^2}, \qquad z \in E.$$

It follows from this and Proposition 2.2 that

$$\mathbb{P}(X \in B) = \mathbf{E}_\mathbb{Q} \mathbb{P}(Y \in B)$$

for every set $B \in \mathcal{B}_\sigma(E')$. In particular, if $\varepsilon > 0$, then

$$\mathbb{P}(\|X\|_{E'} < \varepsilon) = \mathbf{E}_\mathbb{Q} \mathbb{P}(\|Y\|_{E'} < \varepsilon). \tag{2.8}$$

With the definition of the small deviation function (1.2), equation (2.8) may be rewritten as

$$\phi(X, \varepsilon) = -\log(\mathbf{E}_\mathbb{Q} \exp(-\phi(Y, \varepsilon))). \tag{2.9}$$

## 3. Entropy numbers of random operators

Let $u: E \to L_\alpha(S, \sigma)$ be a decomposed operator represented by a certain random operator $v: E \to \ell_2$ as in Proposition 2.2. Our goal is to compare compactness properties of $u$ with those of $v$ and vice versa. We recall Proposition 3.1 from Li and Linde (2004), which is based on an idea from Marcus and Pisier (1984).

**Proposition 3.1.** *There exist universal constants $\rho, \kappa > 0$ such that, for every $m \in \mathbb{N}$,*

$$\mathbb{Q}\left\{ n^{1/\alpha - 1/2} \frac{e_n(v)}{e_n(u)} \geq \rho : n \geq m \right\} \geq 1 - \kappa e^{-m}. \tag{3.1}$$

The proof of Proposition 3.1 rests on the fact that, for each fixed $z \in E$, the non-negative random variable $\|v(z)\|_2^2 / \|u(z)\|_\alpha^2$ is totally skewed $\alpha/2$-stable (in particular, positive). Consequently,

$$\mathbb{Q}\left\{ \frac{\|v(z)\|_2}{\|u(z)\|_\alpha} < \varepsilon \right\} \leq \exp(-c\varepsilon^{-1/(1/\alpha - 1/2)})$$

for some $c > 0$. In order to verify that, similarly to (3.1), an opposite inequality between $e_n(v)$ and $e_n(u)$ holds, this approach does not work. Note that, by the well-known tail behavior of stable random variables, we only get

$$\mathbb{Q}\left\{ \frac{\|v(z)\|_2}{\|u(z)\|_\alpha} > t \right\} \approx t^{-\alpha}$$

as $t \to \infty$. Yet this is far too weak for proving $e_n(v) \leq c n^{-1/\alpha + 1/2} e_n(u)$ on a set of large $\mathbb{Q}$-measure.

Therefore, another approach is needed. In fact, we will prove that the opposite inequality in (3.1) holds (actually, on a set of full $\mathbb{Q}$-measure) if $u$ is replaced by $u_\infty$. Recall that $u$ is assumed to be decomposed by an $E'$-valued function $\varphi$ with $\|\varphi(s)\| \leq 1$ for $s \in S$ and thus $u_\infty$ is well defined.

Before stating and proving this, let us formulate a lemma which is based on the strong law of large numbers. It enables us to replace the random variables $\Gamma_j$ by $j$ in all occurrences where the metric entropy is concerned.

**Lemma 3.2.** *Let $p \in [1, \infty]$. The random diagonal operator $G : \ell_p \to \ell_p$ given by*

$$G: (z_j) \mapsto \left( \left( \frac{\Gamma_j}{j} \right)^{-1/\alpha} z_j \right)$$

*and its inverse are a.s. bounded.*

**Theorem 3.3.** *Let $u$ and $u_\infty$ be as before. We have*

$$\mathbb{Q}\left\{ \limsup_{n \to \infty} n^{1/\alpha - 1/2} \frac{e_{2n-1}(v)}{e_n(u_\infty)} < \infty \right\} = 1.$$

**Proof.** We write $v$ as $D \circ w$, where $w \colon E \to \ell_\infty$ is as in (2.5) and $D \colon \ell_\infty \to \ell_2$ is a diagonal operator as in (2.6). First, note that

$$\|w(z)\|_\infty = \sup_{j \geq 1} |(uz)(V_j)| \leq \|u_\infty(z)\|_\infty$$

for all $z \in E$. Consequently, by Lemma 4.2 in Lifshits and Linde (2002), we get $e_n(w) \leq e_n(u_\infty)$. On the other hand, by Lemma 3.2 above and Theorem 2.2 in Kühn (2005), we obtain, for some random constant $c = c_\delta$,

$$e_n(D \colon \ell_\infty \to \ell_2) \leq c n^{-1/\alpha + 1/2}.$$

Thus, we arrive at

$$e_{2n-1}(v) \leq e_n(u_\infty) \cdot e_n(D) \leq c n^{-1/\alpha + 1/2} e_n(u_\infty)$$

for some *random* constant $c = c_\delta > 0$. This completes the proof of Theorem 3.3. $\square$

## 4. A zero-one law and proof of the main result

### 4.1. A zero-one law for the random operator

In order to proceed in using Theorem 3.3 in the same way as Proposition 3.1 is used in Li and Linde (2004), we have to overcome one essential difficulty. Namely, note that the constant $\rho$ in Proposition 3.1 is not random. Contrary to this, the limit in Theorem 3.3 is a random variable. Our next objective is to show that this random variable is, in fact, almost surely constant.

**Proposition 4.1.** *Let $u$ and $v$ be as in (2.3). For any sequence $(a_n)$ that is regularly varying at infinity, there exists a $C \in [0, \infty]$ such that, $\mathbb{Q}$-almost surely,*

$$\limsup_{n \to \infty} a_n e_n(v) = C.$$

*The same holds for the limit inferior.*

The main idea is to show that $\limsup_{n \to \infty} a_n e_n(v)$ is measurable with respect to the terminal $\sigma$-field and thus a.s. constant. For this purpose, it is sufficient to show that the asymptotic behavior of the entropy of an arbitrary operator $w$ mapping from $\ell_2$ (and thus of the dual $w'$ by Artstein *et al.* (2004)) does not depend on the first components. Due to the special structure of the random operator $v$, the proof of Proposition 4.1 is a direct consequence of the following lemma.

**Lemma 4.2.** *Let $w \colon l_2 \to E$ be some operator and let $P \colon l_2 \to l_2$ be an orthogonal projection of finite rank. Then, for any sequence $(a_n)$ that is regularly varying at infinity,*

$$\limsup_{n \to \infty} a_n e_n(w) = \limsup_{n \to \infty} a_n e_n(w \circ P^\perp).$$

*The same holds for the limit inferior.*

**Proof.** Let $0 < \varepsilon < 1$ and choose a sequence $(k_n)$ of integers such that $k_n \leq n$ and $k_n/n \to \varepsilon$ as $n \to \infty$. It follows that

$$e_n(w) \leq e_{k_n}(w \circ P) + e_{n-k_n+1}(w \circ P^\perp). \tag{4.1}$$

Let us define $m := \text{rank}(P)$. Then, by estimate 1.3.36 in Carl and Stephani (1990),

$$e_{k_n}(w \circ P) \leq c\|w\|2^{-(k_n-1)/m},$$

hence $\lim_{n\to\infty} a_n e_{k_n}(w \circ P) = 0$, for any regularly varying sequence $(a_n)$.

Write $a_n = n^\beta L(n)$ for some $\beta \in \mathbb{R}$ and a slowly varying function $L$. By (4.1), we obtain

$$\limsup_{n\to\infty} a_n e_n(w) \leq \limsup_{n\to\infty} a_n e_{n-k_n+1}(w \circ P^\perp)$$

$$= \limsup_{n\to\infty} \left(\frac{n}{n-k_n+1}\right)^\beta \frac{L(n)}{L(n-k_n+1)} a_{n-k_n+1} e_{n-k_n+1}(w \circ P^\perp)$$

$$\leq \left(\frac{1}{1-\varepsilon}\right)^\beta \cdot 1 \cdot \limsup_{n\to\infty} a_{n-k_n+1} e_{n-k_n+1}(w \circ P^\perp)$$

$$= \left(\frac{1}{1-\varepsilon}\right)^\beta \limsup_{n\to\infty} a_n e_n(w \circ P^\perp).$$

Letting $\varepsilon$ tend to zero, it follows that

$$\limsup_{n\to\infty} a_n e_n(w) \leq \limsup_{n\to\infty} a_n e_n(w \circ P^\perp).$$

In order to see the opposite estimate, start with

$$e_n(w \circ P^\perp) \leq e_{k_n}(w \circ P) + e_{n-k_n+1}(w),$$

hence,

$$e_{n-k_n+1}(w) \geq e_n(w \circ P^\perp) - e_{k_n}(w \circ P)$$

and proceed exactly as before. □

### 4.2. Proof of the main result

**Proof of Theorem 1.3.** The assumption is that $e_n(u_\infty) \preceq n^{1/\alpha - 1/\tau - 1} L(n)$. Consequently, by Theorem 3.3 and Proposition 4.1, there is a finite constant $C \geq 0$ such that, $\mathbb{Q}$-a.s.,

$$\limsup_{n\to\infty} n^{1/\tau + 1/2} L(n)^{-1} e_n(v) = C,$$

which, by Proposition 1.1, implies that, $\mathbb{Q}$-a.s.,

$$\limsup_{\varepsilon \to 0} \varepsilon^\tau L(1/\varepsilon)^{-\tau} \phi(Y,\varepsilon) \leq C'$$

for some $C' < \infty$. Consequently, there exists a constant $C'' \geq 0$ such that, $\mathbb{Q}$-almost surely,

$$\phi(Y,\varepsilon) \leq C'' \varepsilon^{-\tau} L(1/\varepsilon)^\tau \tag{4.2}$$

whenever $\varepsilon < \varepsilon_0$ for some random $\varepsilon_0 > 0$. Thus, we find a non-random $\varepsilon_1 > 0$ such that (4.2) holds for $\varepsilon < \varepsilon_1$ on a set of $\mathbb{Q}$-measure larger than $1/2$. Doing so, it follows that

$$\mathrm{e}^{-\phi(X,\varepsilon)} = \mathbf{E}_\mathbb{Q} \mathrm{e}^{-\phi(Y,\varepsilon)} \geq \tfrac{1}{2} \mathrm{e}^{-C'' \varepsilon^{-\tau} L(1/\varepsilon)^\tau}$$

whenever $\varepsilon < \varepsilon_1$. Hence,

$$\phi(X,\varepsilon) \preceq \varepsilon^{-\tau} L(1/\varepsilon)^\tau,$$

as asserted.　□

**Remark.** The relation between $u$ and $X$ in (1.1) is homogeneous. Thus, Theorem 1.3 can be slightly improved as follows. There exists a constant $c_0 > 0$ depending only on $\alpha$, $\tau$ and $L$ such that, whenever $u$ and $X$ are related via (1.1), it follows that

$$\limsup_{n \to \infty} n^{-1/\alpha + 1/\tau + 1} L(n)^{-1} e_n(u_\infty) =: C$$

implies

$$\limsup_{\varepsilon \to 0} \varepsilon^\tau L(1/\varepsilon)^{-\tau} \phi(X,\varepsilon) \leq c_0 C^\tau.$$

## 5. Open questions

In this paper, we mainly deal with four different objects. The first object is a decomposed operator $u$ from $E$ to $L_\alpha(S,\sigma)$, the second is an $E'$-valued S$\alpha$S random vector $X$ generated by $u$ via (1.1), the third object is the random operator $v = v_\delta$ from $E$ into $\ell_2$ constructed by (2.4) and, finally, we consider the $E'$-valued centered random Gaussian vector $Y = Y_\delta$ associated with $v_\delta$. Several important questions about the relations between these objects remain open.

### 5.1. Question 1

Probably the most interesting set of open questions is whether or not the random operator $v = v_\delta$ (resp., the associated Gaussian vectors $Y_\delta$) determine the small deviation behavior of $X$. In view of (2.8) or (2.9), this is very likely, so we formulate the following conjecture.

**Conjecture 5.1.** *Let $X$ be an $E'$-valued symmetric $\alpha$-stable vector generated by an operator $u : E \to L_\alpha(S, \sigma)$ and denote by $v = v_\delta$ the random operator associated with $u$ via* (2.4). *Let $Y = Y_\delta$ be the corresponding Gaussian vector generated by $v$. Let $\tau > 0$ and let $L$ be a function that is slowly varying at infinity such that $L(t) \approx L(t^p)$ for all $p > 0$. The following equivalences then hold:*

(a) *we have*

$$\phi(Y_\delta, \varepsilon) \succeq \varepsilon^{-\tau} L(1/\varepsilon)^\tau, \quad \mathbb{Q}\text{-}a.s. \quad \Leftrightarrow \quad \phi(X, \varepsilon) \succeq \varepsilon^{-\tau} L(1/\varepsilon)^\tau; \quad (5.1)$$

(b) *we have*

$$\phi(Y_\delta, \varepsilon) \preceq \varepsilon^{-\tau} L(1/\varepsilon)^\tau, \quad \mathbb{Q}\text{-}a.s. \quad \Leftrightarrow \quad \phi(X, \varepsilon) \preceq \varepsilon^{-\tau} L(1/\varepsilon)^\tau. \quad (5.2)$$

*Remark.* Note that, by Proposition 1.1, the left-hand estimate in (5.1) follows from $e_n(v_\delta) \succeq n^{-1/2-1/\tau} L(n)$ a.s. Furthermore, also by Proposition 1.1, observe that the left-hand estimate in (5.2) is equivalent to $e_n(v_\delta) \preceq n^{-1/2-1/\tau} L(n)$ a.s.

Moreover, if we had the regularity of $\phi(Y_\delta, \varepsilon)$ in the sense of Proposition 1.1(a), we could even conclude from the left-hand estimate in (5.1) that $e_n(v_\delta) \succeq n^{-1/2-1/\tau} L(n)$. Also, it might be that the additional regularity condition in Proposition 1.1(a) is not needed *there*.

**Partial proof of Conjecture 5.1.** We can only prove the implications "$\Leftarrow$" in assertion (a) and "$\Rightarrow$" in assertion (b). The two other, more interesting, assertions remain open.

*Proof of the implication "$\Leftarrow$" in* (a). Suppose that

$$\mathbb{P}(\|X\|_{E'} < \varepsilon) \leq \exp(-c\varepsilon^{-\tau} L(1/\varepsilon)^\tau).$$

This means that

$$\mathbf{E}_\mathbb{Q} \mathbb{P}(\|Y_\delta\|_{E'} < \varepsilon) \leq \exp(-c\varepsilon^{-\tau} L(1/\varepsilon)^\tau).$$

Hence, by the Chebyshev inequality for any $c_1 < c$, it follows that

$$\mathbb{Q}\{\delta \in \Delta : \mathbb{P}(\|Y_\delta\|_{E'} < \varepsilon) \geq \exp(-c_1 \varepsilon^{-\tau} L(1/\varepsilon)^\tau)\} \leq \frac{\mathbf{E}_\mathbb{Q} \mathbb{P}(\|Y_\delta\|_{E'} < \varepsilon)}{\exp(-c_1 \varepsilon^{-\tau} L(1/\varepsilon)^\tau)}$$

$$\leq \exp(-(c - c_1)\varepsilon^{-\tau} L(1/\varepsilon)^\tau).$$

Applying the Borel–Cantelli lemma to the sequence $\varepsilon_n = 2^{-n}$, we get that

$$\mathbb{P}(\|Y_\delta\|_{E'} < \varepsilon_n) \leq \exp(-c_1 \varepsilon_n^{-\tau} L(1/\varepsilon_n)^\tau),$$

$\mathbb{Q}$-almost surely, for all $n > n(\delta)$. By properties of regular varying functions, it follows that

$$\mathbb{P}(\|Y_\delta\|_{E'} < \varepsilon) \leq \exp(-c_1 2^{-\tau - 1} \varepsilon^{-\tau} L(1/\varepsilon)^\tau),$$

$\mathbb{Q}$-almost surely, for all $\varepsilon < \varepsilon(\delta)$. Yet this is equivalent to the required estimate $\phi(Y_\delta, \varepsilon) \succeq \varepsilon^{-\tau} L(1/\varepsilon)^\tau$.

*Proof of the implication "$\Rightarrow$" in* (b). This has, in fact, already been done as a step in the proof of Theorem 1.3; see (4.2) and the steps thereafter. □

### 5.2. Question 2

Another interesting question is how the small deviation results depend on the choice of the generating operator $u$. Recall that $u$ is not unique at all. Therefore, the following question is very natural: let $u$ and $\tilde{u}$ be two operators generating the same S$\alpha$S vector $X$, that is, $\|u(z)\|_\alpha = \|\tilde{u}(z)\|_\alpha$ for all $z \in E$. Let $v$ and $\tilde{v}$ be the corresponding random operators. Is it true that

$$e_n(v) \preceq n^{-1/2-1/\tau} L(n), \quad \mathbb{Q}\text{-a.s.} \quad \Leftrightarrow \quad e_n(\tilde{v}) \preceq n^{-1/2-1/\tau} L(n), \quad \mathbb{Q}\text{-a.s.},$$

$$e_n(v) \succeq n^{-1/2-1/\tau} L(n), \quad \mathbb{Q}\text{-a.s.} \quad \Leftrightarrow \quad e_n(\tilde{v}) \succeq n^{-1/2-1/\tau} L(n), \quad \mathbb{Q}\text{-a.s.?}$$

If Conjecture 5.1 holds, then the answer to both questions is affirmative.

When comparing the possible choice of the generating operator $u$, it is worthwhile to note that their variety can be reduced to the following *standard* family. Let $\partial U$ be the unit sphere in $E'$. Recall that for every $S\alpha S$ vector $X$ in $E'$, there exists a unique finite *symmetric* measure $\nu$ concentrated on $\partial U$ such that

$$\mathbf{E} e^{i\langle z, X \rangle} = \exp\left\{ -\int_{E'} |\langle z, x \rangle|^\alpha \, d\nu(x) \right\}, \qquad z \in E. \tag{5.3}$$

The measure $\nu$ is usually called the *spectral measure* of $X$ (see Linde (1986) for further details). Now, let $\tilde{\nu}$ be any measure on $E'$ satisfying the following condition: for any measurable $A \subseteq \partial U$, we have

$$\nu(A) = \frac{1}{2} \int_{\{x \, : \, x/\|x\| \in A\}} \|x\|^\alpha \, d\tilde{\nu}(x) + \frac{1}{2} \int_{\{x \, : \, x/\|x\| \in -A\}} \|x\|^\alpha \, d\tilde{\nu}(x). \tag{5.4}$$

Take $(S, \sigma) = (E', \tilde{\nu})$ and let $u : E \to L_\alpha(E', \tilde{\nu})$ be defined by $(uz)(x) = \langle z, x \rangle$. We then have

$$\|uz\|_\alpha^\alpha = \int_{E'} |\langle z, x \rangle|^\alpha \tilde{\nu}(dx) = \int_{\partial U} |\langle z, x \rangle|^\alpha \, d\nu(x)$$

and the representation condition (1.1) is verified in view of (5.3). We call such representations standard ones.

Obviously, the spectral measure itself satisfies condition (5.4) and provides one possible standard representation. Actually, $\nu$ is the only symmetric measure concentrated on $\partial U$ satisfying (5.3).

Any operator representation can be reduced to a standard one. Indeed, take any representing operator $u : E \to L_\alpha(S, \sigma)$. Let $\varphi : S \to E'$ be a function decomposing $u$, such that (2.1) holds. We then let $\tilde{\nu}$ be the distribution of $\varphi$, namely,

$$\tilde{\nu}(A) = \sigma\{s \in S : \varphi(s) \in A\}, \qquad A \subseteq E'.$$

We claim that the random operators coming from $u$ and from the standard representation associated with $\tilde{\nu}$ have the same distribution and thus possess identical probabilistic properties. Indeed, in the first case, we have

$$v_\delta(z) = ((uz)(V_j)\Gamma_j^{-1/\alpha})_{j=1}^\infty = (\langle z, \varphi(V_j)\rangle \Gamma_j^{-1/\alpha})_{j=1}^\infty,$$

where the $V_j$ are $S$-valued and i.i.d. distributed according to the normalized measure $\sigma$. In the second case, we have

$$\tilde{v}_\delta(z) = (\langle z, \tilde{V}_j\rangle \Gamma_j^{-1/\alpha})_{j=1}^\infty,$$

where $\tilde{V}_j$ are $E'$-valued and i.i.d. distributed according to the normalized measure $\tilde{\nu}$. Clearly, these two sequences are equidistributed.

### 5.3. Question 3

A third question of importance is how compactness properties of $u : E \to L_\alpha(S, \sigma)$ carry over to those of the random operator $v$. To make this more precise, the following definition is useful.

Let $(S, \sigma)$ be a finite measure space and let $0 < \alpha < 2$. If $u : E \to L_\alpha(S, \sigma)$ is decomposed, its *n-th entropy gap* is defined by

$$\mathcal{G}_n(u) := \frac{e_n(v)}{n^{-1/\alpha+1/2} e_n(u)}. \tag{5.5}$$

Here, $v : E \to \ell_2$ is the random operator constructed from $u$ via (2.4). Note that this gap is random.

In view of (3.1), we know that there exists a constant $\rho > 0$ depending only on $\alpha$ such that

$$\mathbb{Q}\Big(\liminf_{n\to\infty} \mathcal{G}_n(u) \geq \rho\Big) = 1.$$

The next result shows that the behavior of the entropy gap is important for our investigations. The assertion follows easily by the methods used to prove Theorem 1.3.

**Proposition 5.2.** *Let $u : E \to L_\alpha(S, \sigma)$ be such that*

$$\mathbb{Q}\Big(\limsup_{n\to\infty} \mathcal{G}_n(u) < \infty\Big) = 1. \tag{5.6}$$

*Under this assumption, the following implication is valid:*

$$e_n(u) \preceq n^{1/\alpha - 1/\tau - 1} L(n) \quad \Rightarrow \quad \phi(X, \varepsilon) \preceq \varepsilon^{-\tau} L(1/\varepsilon)^\tau.$$

Combining Proposition 5.2 with the results from Aurzada (2007a, 2007b), it follows that condition (5.6) cannot be true for arbitrary operators $u$. Thus, the following questions naturally arise:

1. Under what conditions on $u$ is (5.6) satisfied?
2. Given an increasing sequence $(a_n)$ of positive numbers, we say that $\mathcal{G}_n(u)$ has order at most $a_n$ provided that

$$\mathbb{Q}\left(\limsup_{n\to\infty} \frac{\mathcal{G}_n(u)}{a_n} < \infty\right) = 1.$$

One may then ask how large the order of the entropy gap may be. In the next section, we will answer this question for a special class of operators.

## 6. The entropy gap for diagonal operators

### 6.1. Introduction

Recall from Aurzada (2007b) that diagonal operators were used in order to construct the counterexamples mentioned in Proposition 1.2(c). Therefore, special attention should be paid to the investigation of the entropy gap for this type of operators.

First, let us describe which random vector corresponds to a diagonal operator $D$. Namely, let $(\xi_n)$ be a sequence of independent standard S$\alpha$S random variables and let $(\vartheta_n)$ be positive, decreasing coefficients such that $\vartheta_n \to 0$. The random sequence $X = (\vartheta_n \xi_n) \in \ell_{p'}$ is then generated by the diagonal operator $D: \ell_p \to \ell_\alpha$ given by $(z_n) \mapsto (\vartheta_n z_n)$. In order to apply our former results to random vectors of this type, we have to generate them by operators $u$ mapping $\ell_p$ into $L_\alpha(S, \sigma)$ with *finite* measure $\sigma$. This can be done in many different ways. Depending on the special representation, we shall get upper and lower estimates for the entropy gap for certain operators generating vectors $X$ as above. As a consequence, we will see that for these special operators $u$, the entropy gap $\mathcal{G}_n(u)$ defined in (5.5) is (a) not necessarily bounded, but (b) cannot be arbitrarily large. Point (a) will be addressed in Section 6.2 and point (b) in Section 6.3. Actually, we find an integral test exactly describing the possible behavior of $\mathcal{G}_n(u)$ for $u$ related to diagonal operators.

### 6.2. Upper bound for the entropy gap

In this subsection, we work with the following representation of random vectors with values in $\ell_{p'}$ with independent components. Let $S = [0,1]$ and let $\sigma = |\cdot|$ be the Lebesgue measure. We can then define $u: l_p \to L_\alpha[0,1]$ by

$$(z_n) \mapsto \sum_{n=1}^{\infty} \vartheta_n z_n \frac{\mathbb{1}_{A_n}}{|A_n|^{1/\alpha}}, \tag{6.1}$$

where the $A_n$ are disjoint sets in $[0,1]$. It is easy to calculate that this operator generates $X = (\vartheta_n \xi_n)$.

By Theorem 2.2 in Kühn (2005), assuming that $\vartheta_n \approx \vartheta_{2n}$, $|A_n| \approx |A_{2n}|$,

$$\sup_{n \geq k} \left(\frac{n}{k}\right)^a \frac{\vartheta_n}{\vartheta_k} < \infty \qquad \text{for some } a > [1/\alpha - 1/p]_+ \quad \text{and} \tag{6.2}$$

$$\sup_{n \geq k} \left(\frac{n}{k}\right)^b \frac{\vartheta_n |A_n|^{-1/\alpha}}{\vartheta_k |A_k|^{-1/\alpha}} < \infty \qquad \text{for some } b > 0, \tag{6.3}$$

we have

$$e_n(u) \approx \vartheta_n n^{1/\alpha - 1/p} \quad \text{and} \quad e_n(u_\infty) \approx \vartheta_n |A_n|^{-1/\alpha} n^{-1/p}.$$

This yields, by Theorem 3.3, that

$$\mathcal{G}_n(u) = \frac{e_n(v)}{n^{-1/\alpha + 1/2} e_n(u)} \preceq \frac{e_n(u_\infty)}{e_n(u)} \approx \frac{|A_n|^{-1/\alpha}}{n^{1/\alpha}}.$$

Note that we are free in the choice of the sets $A_n$ as long as they are small enough to fit into $[0,1]$. We express this by means of the following integral test.

**Proposition 6.1.** *Let $(d_n)$ and $(\vartheta_n)$ be monotone sequences such that $d_n \approx d_{2n}$, $\vartheta_n \approx \vartheta_{2n}$ and let the regularity conditions (6.2) and (6.3) be valid with $|A_n| := c d_n^{-\alpha} n^{-1}$. Moreover, let*

$$\sum_{n=1}^\infty \frac{d_n^{-\alpha}}{n} < \infty. \tag{6.4}$$

*Then, for the operator $u: \ell_p \to L_\alpha[0,1]$ defined by (6.1), the entropy gap is of order at most $d_n$.*

Conditions (6.2) and (6.3) and the doubling condition are to ensure a certain regularity of the sequence. They are merely technical and due to the application of the results from Kühn (2005). Note that essentially all sequences of interest satisfy these conditions.

Let us illustrate with an example how Proposition 6.1 works.

*Example.* Choose $d_n := (\log n)^{\gamma/\alpha}$ for some $\gamma > 1$. Clearly, (6.3) and (6.4) hold in this case, as well as the doubling condition. Consequently, for all sequences $\vartheta_n$ satisfying (6.2) and $\vartheta_n \approx \vartheta_{2n}$, we have

$$\mathcal{G}_n(u) \preceq (\log n)^{\gamma/\alpha},$$

where $u$ is defined by (6.1). Note that this is valid for any $\gamma > 1$. Of course, we may also take $d_n := (\log n)^{1/\alpha} (\log \log n)^{\gamma/\alpha}$ for some $\gamma > 1$ or another regular sequence $(d_n)$ satisfying condition (6.4). Any such summable sequence that is sufficiently regular yields an upper bound for $\mathcal{G}_n(u)$.

### 6.3. The entropy gap of embedding operators

In this subsection, we give examples in which the entropy gap is unbounded. More precisely, we show that it can increase as at least $(\log n)^{1/\alpha}$ and even slightly faster for operators generating stable vectors in $\ell_{p'}$ with independent components.

Here, we use another representation for the generating operator of $X$. Namely, we will use $S = \mathbb{N}$ and the measure $\sigma$ is given by the weights $\sigma_n =: \sigma(\{n\})$, where $\sigma_n = \vartheta_n^\alpha$ in the notation of Section 6.2. We then consider the embedding operator $u : \ell_p \to L_\alpha(\mathbb{N}, \sigma)$. It is straightforward to check that $u$ generates the $\ell_{p'}$-valued random vector $X = (\sigma_n^{1/\alpha} \xi_n)$.

Let us first look at the representation (2.7) and the random operator $v$ occurring there. For this purpose, consider the random operator $w : \ell_p \to \ell_\infty$ defined as follows. Let $(\mathbf{e}_k)$ be the standard basis in $\ell_p$ and let $(V_j)$ be i.i.d. $\mathbb{N}$-valued random variables distributed according to $\sigma$. We set

$$w(\mathbf{e}_k) = \sum_{j : V_j = k} \mathbf{e}_j, \quad \text{that is,} \quad w(z) = \sum_k \sum_{j : V_j = k} z_k \mathbf{e}_j.$$

One can interpret this object as a random partitioning of $\mathbb{N}$ into sets $B_k = \{j : V_j = k\}$. Every point is put into $B_k$ independently of other points and with the same probability $\sigma_k$ for all points.

We combine this operator $w$ with the diagonal operator $D : \ell_\infty \to \ell_2$ possessing the diagonal $j^{-1/\alpha}$. The result of the combination is the operator $D \circ w : \ell_p \to \ell_2$ acting as

$$(D \circ w)(z) = \sum_k \sum_{j : V_j = k} z_k j^{-1/\alpha} \mathbf{e}_j.$$

Since

$$\|(D \circ w)(z)\|_2^2 = \sum_k z_k^2 \sum_{j : V_j = k} j^{-2/\alpha},$$

the operator $D \circ w$ is clearly isomorphic (as regards its image and hence its compactness properties) to a diagonal operator with the random diagonal

$$\lambda_k = \left( \sum_{j : V_j = k} j^{-2/\alpha} \right)^{1/2}. \tag{6.5}$$

Once it comes to the entropy numbers, by Lemma 3.2, we can replace the operator $\tilde{D} : \ell_\infty \to \ell_2$ with diagonal $\Gamma_j^{-1/\alpha}$, where the $\Gamma_j$ are as above, by the operator $D$.

Finally, we note that $\tilde{D} \circ w = v$ with $v$ from the mixture (2.7). Recall that we are interested in the relation between $e_n(v)$ and $e_n(u)$. From the above arguments, it is clear that we can also consider $e_n(D \circ w)$ instead of $e_n(v)$.

For this purpose, we are interested in the decreasing rearrangement $\lambda_k^*$ of the sequence in (6.5). Before giving a sharp and precise result, let us illustrate the situation by means of two basic cases:

(a) $\sigma_k = ck^{-1}(\log k)^{-\nu}$ with $\nu > 1$;
(b) $\sigma_k = ck^{-a}(\log k)^{-\nu}$ with $a > 1, \nu \in \mathbb{R}$.

**Proposition 6.2.** *For the above cases, we obtain the following results.*

- *In case* (a), *we have* $\lambda_k^* \approx k^{-1/\alpha}(\log k)^{-(\nu-1)/\alpha}$ *almost surely. Hence, it follows that* $e_n(D \circ w) \approx n^{-(1/\alpha - 1/2 + 1/p)}(\log n)^{-(\nu-1)/\alpha}$.
- *In case* (b), *we have* $\lambda_k^* \approx k^{-a/\alpha}(\log k)^{-\nu/\alpha}$ *almost surely. Hence, it follows that* $e_n(D \circ w) \approx n^{-(a/\alpha - 1/2 + 1/p)}(\log n)^{-\nu/\alpha}$.

We will give the proof of this result after that of Proposition 6.4.

Note that Proposition 6.2 has the following consequence for the entropy gap.

**Corollary 6.3.** *We have*

$$\mathcal{G}_n(u) \approx \frac{e_n(D \circ w)}{n^{-1/\alpha + 1/2} e_n(u)} \approx \begin{cases} (\log n)^{1/\alpha}, & \text{in case } (a), \\ 1, & \text{in case } (b). \end{cases}$$

This is why one can call case (b) "regular" and case (a) "exceptional".

Our main result concerning possible entropy gaps for embedding operators from $\ell_p$ into $L_\alpha(\mathbb{N}, \sigma)$ is as follows.

**Proposition 6.4.** *Let* $(d_k)$ *be an increasing positive sequence such that* $d_n \approx d_{2n}$,

$$\sum_{k=1}^{\infty} \frac{d_k^{-\alpha}}{k} = \infty \tag{6.6}$$

*and*

$$d_k \succeq (\log k)^{1/2\alpha}. \tag{6.7}$$

*Then, for any* $p \in [1, \infty)$, *there exists a probability measure* $\sigma$ *on* $\mathbb{N}$ *such that for the embedding* $u: \ell_p \to L_\alpha(\mathbb{N}, \sigma)$, *we have*

$$\mathcal{G}_n(u) \succeq d_n.$$

The technical assumption (6.7) could probably be avoided, but it is not an obstacle for considering interesting examples of gaps which are of order $(\log n)^{1/\alpha}$ and larger.

Note that the integral test in (6.6) is the same as in (6.4) of the previous subsection. Therefore, both results are sharp.

**Proof of Proposition 6.4.** Suppose that some probability measure $\sigma$ on $\mathbb{N}$ is given and let the $V_j$ be independent, $\sigma$-distributed integers. As before, $\sigma_k := \sigma(\{k\})$. For any integer $m$, consider the set $\{V_j, j \leq m\}$. The key question we address now is how many *different* values there are in this set for large $m$. Consider the random events $G_k = \{\exists j \leq m : V_j = k\}$. Let $I_k = \mathbb{1}_{G_k}$. Clearly,

$$\mathbb{P}(G_k) = 1 - \mathbb{P}(\bar{G}_k) = 1 - (1 - \sigma_k)^m \sim m\sigma_k \tag{6.8}$$

whenever $m\sigma_k \to 0$.

We are interested in the behavior of the number of different values

$$N_m := \sum_{k=1}^{\infty} I_k. \tag{6.9}$$

Let us look at the variance of $N_m$. Note that the $G_k$ are negatively dependent: for $k_1 \neq k_2$, we have $\mathbb{P}(G_{k_1}|G_{k_2}) \leq \mathbb{P}(G_{k_1})$. In other words, $\mathbb{P}(G_{k_1} \cap G_{k_2}) \leq \mathbb{P}(G_{k_1})\mathbb{P}(G_{k_1})$. The latter relation can be also written as $\text{cov}(I_{k_1}, I_{k_2}) \leq 0$. It follows that

$$\text{Var}(N_m) \leq \sum_{k=1}^{\infty} \text{cov}(I_k, I_k) \leq \sum_{k=1}^{\infty} \mathbb{P}(G_k) = \mathbf{E}N_m.$$

If we define $a_k := \frac{d_k^{-\alpha}}{k}$, then $A_n := \sum_{k=1}^{n} a_k \nearrow \infty$. We set $\sigma_k = ca_{k+1} \exp(-A_k)$, where $c$ is a normalizing constant such that $\sum_{k=1}^{\infty} \sigma_k = 1$. Note that since $a_k \to 0$,

$$a_{k+1} \exp(-A_k) \sim (1 - \exp(-a_{k+1})) \exp(-A_k) = \exp(-A_k) - \exp(-A_{k+1})$$

really forms a convergent series, thus the choice of the normalizer $c$ is possible.

By construction, $\sigma_k$ is a decreasing sequence and, as before, the operator $u$ is defined as embedding from $\ell_p$ into $L_\alpha(\mathbb{N}, \sigma)$. It corresponds to the diagonal operator from $\ell_p$ to $\ell_\alpha$ with diagonal $(\vartheta_n)$, where $\vartheta_n = \sigma_n^{1/\alpha}$, as usual.

Now, define the tails

$$T_n := \sum_{k=n}^{\infty} \sigma_k \sim c \sum_{k=n}^{\infty} (\exp(-A_k) - \exp(-A_{k+1})) = c \exp(-A_n). \tag{6.10}$$

Hence,

$$\frac{\sigma_n}{T_n} \sim a_{n+1}. \tag{6.11}$$

We will now use the numbers $N_m$, as defined in (6.9). Since $a_k = O(1/k)$ and $A_k \to \infty$, our construction yields

$$\max_{k \geq m} \sigma_k = \sigma_m = o(1/m).$$

It follows from (6.8) that

$$\mathbf{E}N_m \geq \sum_{k=m}^{\infty} \mathbb{P}(G_k) \sim m \sum_{k=m}^{\infty} \sigma_k = mT_m$$

and we obtain $\mathbf{E}N_m \geq (1-\varepsilon)mT_m$ for any $\varepsilon > 0$ and all large $m$. Also, recall that $\mathrm{Var}(N_m) \leq \mathbf{E}N_m$, hence, by the Chebyshev inequality,

$$\mathbb{P}(N_m \leq (1-2\varepsilon)mT_m) \leq \mathbb{P}(N_m \leq (1-\varepsilon)\mathbf{E}N_m) \leq \frac{1}{\varepsilon^2(1-\varepsilon)mT_m}.$$

We now show that $T_m$ decreases rather slowly, leading to a convergent series in the Borel–Cantelli lemma. Indeed, since $a_k = \mathrm{o}(1/k)$, we have $A_n = \mathrm{o}(\log n)$. Hence, for large $n$, by (6.10),

$$T_n \sim c\exp(-A_n) \succeq \exp(-\log n/2) = n^{-1/2}.$$

By the Borel–Cantelli lemma, it follows that $N_m \geq (1-2\varepsilon)mT_m$ for large $m$ along exponential sequences. Using the fact that $N_m$ and $T_m$ are monotone, we easily obtain that

$$N_m \geq (1-3\varepsilon)mT_m \tag{6.12}$$

for all large $m$.

Recall that if $j, k, m$ are such positive integers such that $j \leq m$ and $V_j = k$, then $\lambda_k \geq j^{-1/\alpha} \geq m^{-1/\alpha}$. Therefore, (6.12) can be written as

$$\lambda^*_{\lceil (1-2\varepsilon)mT_m \rceil} \geq m^{-1/\alpha}. \tag{6.13}$$

Let $m = m(n) = \lceil n/T_n \rceil$ and observe that for some $\eta > 0$,

$$\liminf_{n \to \infty} \frac{T_m}{T_n} > \eta. \tag{6.14}$$

There then exists a small $\varepsilon > 0$ such that for all large $n$,

$$\eta n \leq (1-2\varepsilon)\frac{T_m}{T_n}n \leq (1-2\varepsilon)mT_m.$$

By applying (6.13), we get

$$\lambda^*_{\lceil \eta n \rceil} \geq \lambda^*_{\lceil (1-2\varepsilon)mT_m \rceil} \succeq (n/T_n)^{-1/\alpha}.$$

Next, note that by (6.10),

$$1 \leq \frac{T_{\lceil \eta n \rceil}}{T_n} \sim \exp\left(\sum_{k=\lceil \eta n \rceil+1}^{n} a_k\right) \leq \exp\left(d^{-\alpha}_{\lceil \eta n \rceil} \sum_{k=\lceil \eta n \rceil+1}^{n} \frac{1}{k}\right) \to 1. \tag{6.15}$$

It follows that

$$\lambda^*_{\lceil \eta n \rceil} \succeq (\lceil \eta n \rceil/T_{\lceil \eta n \rceil})^{-1/\alpha}.$$

By changing the notation and using (6.11), we obtain

$$\lambda^*_n \succeq (T_n/n)^{1/\alpha} \approx \left(\frac{\sigma_n}{na_n}\right)^{1/\alpha} = \sigma_n^{1/\alpha}d_n. \tag{6.16}$$

We continue by proving (6.14), which, by (6.10), is equivalent to

$$\limsup_{n \to \infty} \sum_{k=n+1}^{m} a_k < \infty.$$

However, assumption (6.7) yields that

$$\sum_{k=n+1}^{m} a_k \preceq \sum_{n+1}^{m} \frac{1}{(\log k)^{1/2} k} \preceq (\log m)^{1/2} - (\log n)^{1/2} \frac{|\log T_n|}{(\log n)^{1/2}} \sim \frac{A_n}{(\log n)^{1/2}} \preceq 1.$$

We can now complete our evaluation of the entropy gap by using information about our diagonal operators. First, consider the non-random operator $u$. Recall that $u$ is a diagonal operator with diagonal $(\vartheta_n)$ defined by

$$\vartheta_n = \sigma_n^{1/\alpha} = a_{n+1}^{1/\alpha} \exp(-A_n/\alpha) = (n+1)^{-1/\alpha} d_{n+1}^{-1} \exp(-A_n/\alpha).$$

Note that the second and third factors are decreasing sequences. We see that the standard regularity condition that is necessary to get the entropy behavior is verified, namely, for $a = 1/\alpha$, we have

$$\sup_{n \geq k} \left(\frac{n}{k}\right)^a \frac{\vartheta_n}{\vartheta_k} < \infty.$$

Recall that $(d_n)$ satisfies the doubling condition $d_n \approx d_{2n}$. Moreover, the sequence $\exp(-A_n/\alpha)$ is slowly varying; see (6.15). Hence, $\vartheta_n$ satisfies the doubling condition $\vartheta_n \approx \vartheta_{2n}$. By the aforementioned Theorem 2.2 in Kühn (2005), it follows that

$$e_n(u) \approx \vartheta_n n^{1/\alpha - 1/p'} = \sigma_n^{1/\alpha} n^{1/\alpha - 1/p'}.$$

The same arguments apply to the lower bound (6.16) which we obtained for $\lambda_n^*$. Hence,

$$e_n(v) \succeq \sigma_n^{1/\alpha} d_n n^{1/2 - 1/p'}.$$

It follows that $\mathcal{G}_n(u) \succeq d_n$, as required. □

**Proof of Proposition 6.2.** *Lower bound.* Setting $d_k = (\log k)^{1/\alpha}$ in the previous construction yields the lower bounds for $\lambda_n^*$ in case (a). In case (b), the calculation is quite similar: by direct calculation of the mean, followed by a Borel–Cantelli argument, we get

$$N_m \approx \mathbf{E} N_m \geq c \frac{m^{1/a}}{(\log m)^{\nu/a}}.$$

Since

$$\{N_m > k\} \subseteq \{m^{-1/\alpha} \leq \lambda_{k+1}^*\},$$

we get $\lambda_k^* \geq c k^{-a/\alpha} (\log k)^{-\nu/\alpha}$, as required.

*Upper bound.* We first treat case (a). Let us introduce some notation. Let $r$ be a small number. Let $k_* = \frac{r^{-\alpha/2}}{|\log r|^{\nu-1}}$,

$$F_1(r) = \#\left\{k : \sum_{j\,:\,V_j=k} j^{-2/\alpha} \geq r\right\}$$

and

$$F_2(r) = \#\left\{k \geq k_* : \sum_{j\,:\,j \geq r^{-\alpha/2},\, V_j=k} j^{-2/\alpha} \geq r\right\}.$$

Clearly,

$$F_1(r) \leq F_2(r) + k_* + N(r^{-\alpha/2}), \tag{6.17}$$

where $N(m) = N_m$ was defined in (6.9).

To evaluate $F_2(r)$, we need

$$T_k(r) = \sum_{j\,:\,j \geq r^{-\alpha/2},\, V_j=k} j^{-2/\alpha}$$

and

$$\mathcal{S}(r) = \sum_{k \geq k_*} T_k(r) = \sum_{j\,:\,j \geq r^{-\alpha/2},\, V_j \geq k_*} j^{-2/\alpha}.$$

Note that the latter expression is a weighted sum of independent Bernoulli variables whose parameters are

$$\mathbb{P}(V_j \geq k_*) = \sum_{k \geq k_*} \sigma_k \approx |\log r|^{-(\nu-1)}.$$

We now evaluate the expectation and variance of $\mathcal{S}(r)$. Indeed,

$$\mathbf{E}\mathcal{S}(r) = \sum_{j \geq r^{-\alpha/2}} j^{-2/\alpha} \mathbb{P}(V_j \geq k_*) \approx r^{1-\alpha/2} |\log r|^{-(\nu-1)},$$

$$\operatorname{Var}\mathcal{S}(r) = \sum_{j \geq r^{-\alpha/2}} j^{-4/\alpha} \operatorname{Var} \mathbb{1}_{\{V_j \geq k_*\}} \preceq r^{2-\alpha/2} |\log r|^{-(\nu-1)}.$$

By the Chebyshev inequality,

$$\mathbb{P}(\mathcal{S}(r) \geq 2\mathbf{E}\mathcal{S}(r)) \leq c r^{\alpha/2} |\log r|^{\nu-1}.$$

Again using the asymptotics of $\mathbf{E}\mathcal{S}(r)$ and the trivial inequality $F_2(r) r \leq \mathcal{S}(r)$, we get

$$\mathbb{P}(F_2(r) \geq c_1 r^{-\alpha/2} |\log r|^{-(\nu-1)}) \leq c_2 r^{\alpha/2} |\log r|^{\nu-1}.$$

By the Borel–Cantelli lemma, we conclude that

$$F_2(r) = \mathrm{O}(r^{-\alpha/2}|\log r|^{-(\nu-1)})$$

almost surely, at least along the subsequence $r = 2^{-i}, i = 1, 2, \ldots$. Next, we pass from $F_2$ to $F_1$. To this end, the quantity $N(\cdot)$ in (6.17) should be evaluated. By using (6.8), one easily finds that $\mathbf{E} N_m \sim cm(\log m)^{-(\nu-1)}$ in case (a). Moreover, since $\operatorname{Var} N_m \leq \mathbf{E} N_m$, a Borel–Cantelli argument shows that $N_m \leq 2cm(\log m)^{-(\nu-1)}$ for all sufficiently large $m$. In particular,

$$N(r^{-\alpha/2}) = \mathrm{O}(r^{-\alpha/2}|\log r|^{-(\nu-1)}).$$

It now follows from the definition of $k_*$ and (6.17) that

$$F_1(r) = \mathrm{O}(r^{-\alpha/2}|\log r|^{-(\nu-1)})$$

almost surely along the aforementioned subsequence. However, since $F_1(\cdot)$ is a decreasing function, the statement is also true along $r \to 0$. This means that

$$\#\{k : \lambda_k \geq r\} = \mathrm{O}(r^{-\alpha}|\log r|^{-(\nu-1)}), \qquad r \to 0,$$

which is equivalent to the required estimate

$$\lambda_k^* \leq ck^{-1/\alpha}|\log k|^{-(\nu-1)/\alpha}.$$

Therefore, we are finished with the upper estimate in case (a). For case (b), set $k_* = \frac{r^{-\alpha/2a}}{|\log r|^{\nu/a}}$. By repeating the previous calculations, we subsequently get

$$\mathbf{E}\mathcal{S}(r) \approx r^{1-\alpha/2a}|\log r|^{-\nu/a}, \qquad F_2(r) = \mathrm{O}(r^{-\alpha/2a}|\log r|^{-\nu/a})$$

and

$$N_m = \mathrm{O}(m^{1/a}|\log m|^{-\nu/a}), \qquad F_1(r) = \mathrm{O}(r^{-\alpha/2a}|\log r|^{-\nu/a}),$$

which gives

$$\#\{k : \lambda_k \geq r\} = \mathrm{O}(r^{-\alpha/a}|\log r|^{-\nu/a}), \qquad r \to 0,$$

or $\lambda_k^* \leq ck^{-a/\alpha}(\log k)^{-\nu/\alpha}$, as required. $\square$

# 7. Examples and applications

## 7.1. Application to symmetric $\alpha$-stable processes

### 7.1.1. Symmetric $\alpha$-stable processes

A stochastic process $X = (X(t))_{t \in T}$ indexed by a non-empty set $T$ is said to be S$\alpha$S for some $\alpha \in (0, 2]$ if, for all $t_1, \ldots, t_n \in T$ and all real numbers $\lambda_1, \ldots, \lambda_n$, the real random variable $\sum_{j=1}^n \lambda_j X(t_j)$ is S$\alpha$S-distributed.

We shall restrict ourselves to S$\alpha$S processes possessing an integral representation in the sense of Chapter 13 of Samorodnitsky and Taqqu (1994). We note that all natural examples of S$\alpha$S processes fit into this framework.

In other words, we investigate S$\alpha$S processes $X$ for which there exist a measure space $(S, \sigma)$ and a kernel $K : T \times S \to \mathbb{R}$ such that for each $t \in T$, the function $s \mapsto K(t, s)$ is measurable with

$$\int_S |K(t,s)|^\alpha \, \mathrm{d}\sigma(s) < \infty$$

and for all $\lambda_1, \ldots, \lambda_n \in \mathbb{R}$ and all $t_1, \ldots, t_n \in T$, we have

$$\mathbf{E} \exp\left( \mathrm{i} \sum_{j=1}^n \lambda_j X(t_j) \right) = \exp\left( -\int_S \left| \sum_{j=1}^n \lambda_j K(t_j, s) \right|^\alpha \mathrm{d}\sigma(s) \right). \tag{7.1}$$

Usually, one writes

$$X(t) = \int_S K(t,s) \, \mathrm{d}M(s), \qquad t \in T, \tag{7.2}$$

where $M$ denotes an independently scattered S$\alpha$S random measure with control measure $\sigma$. If $S \subseteq \mathbb{R}$ and $\sigma$ is the Lebesgue measure on $S$, then

$$X(t) = \int_S K(t,s) \, \mathrm{d}Z_\alpha(s), \qquad t \in T.$$

We refer to Samorodnitsky and Taqqu (1994) for more information about integral representations of S$\alpha$S processes.

Now, suppose that $(T, d)$ is a separable metric space endowed with the Borel $\sigma$-field. If the kernel $K$ on $T \times S$ is measurable with respect to the product $\sigma$-field, then $X$ possesses a measurable version. Let $\mu$ be some finite Borel measure on $T$ and suppose that

$$\mathbb{P}(\|X\|_{L_q(T,\mu)} < \infty) = 1 \tag{7.3}$$

for a certain $q \in [1, \infty]$. Recall from Section 11.3 in Samorodnitsky and Taqqu (1994) that, if $q < \infty$, there is a simple condition in terms of the kernel $K$ to verify (7.3). We now regard $X$ as an S$\alpha$S random vector in $L_q(T, \mu)$ and define $p \in [1, \infty]$ by $p' = q$. Then, as proven in Li and Linde (2004), Proposition 5.1, the operator $u: L_p(T, \mu) \to L_\alpha(S, \sigma)$ with

$$(uf)(s) := \int_T K(t,s) f(t) \, d\mu(t), \qquad s \in S, \tag{7.4}$$

generates $X$ in the sense of (1.1). Consequently, by Theorem 1.3, any upper entropy estimate for $u: L_p(T, \mu) \to L_\infty(S, \sigma)$ implies an upper estimate for

$$\phi(X, \varepsilon) = -\log \mathbb{P}(\|X\|_{L_q(T, \mu)} < \varepsilon),$$

with $q = p'$. We summarize these observations as follows.

**Proposition 7.1.** *Let $p \in [1, \infty]$, $K$ be as above and $u$ be as in (7.4). Fix $\tau > 0$ and a slowly varying function $L$ as in Theorem 1.3. Assume that*

$$e_n(u: L_p(T, \mu) \to L_\infty(S, \sigma)) \preceq n^{1/\alpha - 1/\tau - 1} L(n).$$

*Set $q = p'$. Then,*

$$-\log \mathbb{P}(\|X\|_{L_q(T, \mu)} < \varepsilon) \preceq \varepsilon^{-\tau} L(1/\varepsilon)^\tau.$$

Let us illustrate this by means of several concrete examples.

*7.1.2. Hölder operators*

We begin our investigation of the small deviations of S$\alpha$S processes with a quite general approach. To this end, suppose that $(S, d)$ is a compact metric space and let $C(S)$ be the Banach space of (real-valued) continuous functions on $S$. An operator $u: E \to C(S)$ is said to be $\beta$-*Hölder* for some $\beta \in (0, 1]$ provided there exists a constant $c > 0$ such that for all $z \in E$ and all $s_1, s_2 \in S$, it follows that

$$|(uz)(s_1) - (uz)(s_2)| \leq c \|z\|_E d(s_1, s_2)^\beta.$$

Furthermore, let $\varepsilon_n(S)$ be the sequence of covering numbers of $S$ (with respect to the metric $d$). The basic result about compactness properties of Hölder operators is as follows (see Carl and Stephani (1990)).

**Proposition 7.2.** *Let $\mathcal{H}$ be a Hilbert space and let $u: \mathcal{H} \to C(S)$ be a $\beta$-Hölder operator. If $\varepsilon_n(S) \leq h(n)$ for some regularly varying decreasing function $h$, then we have*

$$e_n(u) \leq c n^{-1/2} h(n)^\beta.$$

We apply this result in our setup. To this end, let $(S, d)$ be as before and suppose that $\sigma$ is a finite Borel measure on $S$.

**Proposition 7.3.** *Let $X$ be an S$\alpha$S vector with values in a Hilbert space $\mathcal{H}$ and let $u\colon\mathcal{H}\to L_\alpha(S,\sigma)$ be an operator generating $X$. Suppose that $u$ is $\beta$-Hölder for some $\beta\in(0,1]$ and that $\varepsilon_n(S)\preceq n^{-\gamma}L(n)$ for some $\gamma>0$ and some slowly varying function $L$, as before. Set $1/\tau=1/\alpha-1/2+\gamma\beta$. We then have*

$$\phi(X,\varepsilon)\preceq \varepsilon^{-\tau}L(1/\varepsilon)^{\beta\tau}.$$

**Proof.** Due to the assumptions, the operator $u_\infty$ maps $\mathcal{H}$ into $C(S)$ and, moreover, it is $\beta$-Hölder. Consequently, Proposition 7.2 applies to $u_\infty$ and yields

$$e_n(u_\infty)\preceq n^{-1/2-\gamma\beta}L(n)^\beta$$

since $\varepsilon_n(S)\preceq n^{-\gamma}L(n)$. From this, the assertion follows immediately from Theorem 1.3. $\square$

Let us apply the preceding result to S$\alpha$S processes with integral representations, as in Section 7.1.1. To this end, suppose that $T$ is a separable metric space with finite Borel measure $\mu$. Assume that an S$\alpha$S process $X=(X(t))_{t\in T}$ has a.s. paths in $L_2(T,\mu)$. If $X$ admits the representation (7.2) with respect to the control measure $\sigma$ on the compact space $(S,d)$, then the process is generated by $u\colon L_2(T,\mu)\to L_\alpha(S,\sigma)$ defined in (7.4). Note that this $u$ is $\beta$-Hölder if and only if there exists some $c>0$ such that, for all $s_1,s_2\in S$, it follows that

$$\left(\int_T |K(t,s_1)-K(t,s_2)|^2\,\mathrm{d}\mu(t)\right)^{1/2}\le cd(s_1,s_2)^\beta. \tag{7.5}$$

Rewriting Proposition 7.3 in this framework, we obtain the following result.

**Corollary 7.4.** *Let*

$$X(t)=\int_S K(t,s)\,\mathrm{d}M(s),\qquad t\in T,$$

*where $M$ has control measure $\sigma$ and $T$ and $S$ are as before. Suppose that (7.5) holds for some $\beta\in(0,1]$. If $\varepsilon_n(S)\preceq n^{-\gamma}L(n)$, then this implies that*

$$-\log\mathbb{P}(\|X\|_{L_2(T,\mu)}<\varepsilon)\preceq \varepsilon^{-\tau}L(1/\varepsilon)^{\beta\gamma},$$

*where $1/\tau=1/\alpha-1/2+\gamma\beta$ as in Proposition 7.3.*

We will show in the next subsection that this leads to sharp estimates in several examples.

### 7.1.3. Riemann–Liouville processes

The symmetric $\alpha$-stable Riemann–Liouville process on $[0,1]$ with Hurst index $H>0$ is usually defined by

$$R_H^\alpha(t):=\int_0^t (t-s)^{H-1/\alpha}\,\mathrm{d}Z_\alpha(s),\qquad 0\le t\le 1,$$

where $Z_\alpha$ is, as above, the symmetric $\alpha$-stable Lévy motion. In other words, $R_H^\alpha$ is the process satisfying (7.1) with $K(t,s) = (t-s)^{H-1/\alpha} \mathbb{1}_{[0,t]}(s)$. The underlying measure space is $T = [0,1]$ endowed with the Lebesgue measure.

The small deviation behavior of $(R_H^\alpha(t))_{0 \le t \le 1}$ was investigated thoroughly in Lifshits and Simon (2005). It was shown that for $H \ge 1/\alpha$, that is, when the process $R_H^\alpha$ is a.s. bounded, for any $q \in [1, \infty]$, it is true that

$$-\log \mathbb{P}(\|R_H^\alpha\|_{L_q[0,1]} < \varepsilon) \sim c\varepsilon^{-1/H} \tag{7.6}$$

for some finite, positive $c = c(H, \alpha)$. One should also mention that for the case $H = 1/\alpha$, the process $R_H^\alpha$ is just a symmetric $\alpha$-stable Lévy motion and its small deviations were studied some time ago in works by Borovkov and Mogul'skiĭ (1991) and Mogul'skiĭ (1974).

In the case $[1/\alpha - 1/q]_+ < H \le 1/\alpha$, the paths are no longer bounded, yet belong to $L_q[0,1]$. For these $H$, the behavior of (7.6) was stated as an open question (see Lifshits and Simon (2005), Section 6.4). Using Theorem 1.3, we can answer this question as follows.

**Proposition 7.5.** *Suppose that $H > [1/\alpha - 1/q]_+$. It then follows that*

$$-\log \mathbb{P}(\|R_H^\alpha\|_{L_q[0,1]} < \varepsilon) \sim c\varepsilon^{-1/H}$$

*for some finite positive $c = c(H, \alpha)$.*

**Proof.** The existence of the positive (but possibly infinite) limit

$$c = \lim_{\varepsilon \to 0}[-\varepsilon^{1/H} \log \mathbb{P}(\|R_H^\alpha\|_{L_q[0,1]} < \varepsilon)]$$

follows from Theorem 4 in Lifshits and Simon (2005). It thus remains to verify the upper estimate. To this end, we set $p := q'$ and consider the integral operator $u \colon L_p[0,1] \to L_\alpha[0,1]$ generating $R_H^\alpha$. This is given by

$$(uf)(s) := \int_s^1 (t-s)^{H-1/\alpha} f(t)\, dt, \qquad 0 \le s \le 1.$$

Since $H > 1/\alpha - 1/q = 1/\alpha - 1/p'$, the operator $u$ even maps $L_p[0,1]$ into $L_\infty[0,1]$. Moreover, introducing the changes of variables $t \mapsto 1-t$ and $s \mapsto 1-s$, it follows that $e_n(u_\infty) = e_n(R^H)$, where $R^H \colon L_p[0,1] \to L_\infty[0,1]$ is (up to a constant) the usual Riemann–Liouville integration operator defined by

$$(R^H f)(s) := \int_0^s (s-t)^{H-1/\alpha} f(t)\, dt. \tag{7.7}$$

As shown in Aurzada and Simon (2007), proof of Lemma 3.9, or in Li and Linde (1999), we have

$$e_n(u_\infty) = e_n(R^H \colon L_p[0,1] \to L_\infty[0,1]) \approx n^{-H-1+1/\alpha}.$$

Thus, the finiteness of $c$ follows from Theorem 1.3. □

**Remark.** The same result holds for linear stable fractional motion

$$X_H^\alpha(t) := \int_{-\infty}^t [(t-s)^{H-1/\alpha} - (-s)_+^{H-1/\alpha}] \, dZ_\alpha(s), \qquad 0 \le t \le 1, 0 < H < 1,$$

since it was shown in Lifshits and Simon (2005) that the difference between $R_H^\alpha$ and $X_H^\alpha$ is irrelevant as far as small deviations are concerned.

**Remark.** The results of Aurzada and Simon (2007) can also be improved to a larger range of $H$ when considering $L_q$-norms, $q < \infty$, by the use of Proposition 7.5.

Finally, let us look at the Riemann–Liouville processes from the point of view of Hölder operators. If $1/\alpha - 1/2 < H \le 1/\alpha + 1/2$, then, as can be easily seen, the corresponding operator $u$ given in (7.7) satisfies (7.5) with $\beta = H - 1/\alpha + 1/2$. Since $S = [0,1]$, we have $\varepsilon_n(S) \approx n^{-1}$. Hence, in that case,

$$1/\tau = 1/\alpha - 1/2 + \beta = H,$$

giving $\tau = 1/H$. This leads to the sharp estimate (compare with Proposition 7.5)

$$-\log \mathbb{P}\left(\int_0^1 |R_H^\alpha(t)|^2 \, dt < \varepsilon^2\right) \preceq \varepsilon^{-1/H}, \qquad (7.8)$$

at least if $1/\alpha - 1/2 < H \le 1/\alpha + 1/2$.

### 7.1.4. Weighted $\alpha$-stable Lévy motion

Let $\rho: [0,1] \to [0, \infty)$ be some (measurable) weight function. The weighted symmetric $\alpha$-stable Lévy motion is defined by

$$X_\rho(t) := \rho(t) Z_\alpha(t), \qquad 0 \le t \le 1.$$

If $1 \le q < \infty$, then there exists a complete characterization of weights such that $X_\rho \in L_q[0,1]$ a.s. (see Example 2 in Li and Linde (2004)).

Using Theorem 1.3, we get the following result.

**Proposition 7.6.** *Let $q \in [1, \infty)$ and assume that $\rho \in L_q[0,1]$. Let $r > 0$ be defined by $1/r := 1/q + 1/\alpha$. There then exist $c_1, c_2 > 0$ such that*

$$\begin{aligned} c_1 \|\rho\|_r &\le \liminf_{\varepsilon \to 0} \varepsilon^\alpha [-\log \mathbb{P}(\|X_\rho\|_q < \varepsilon)] \\ &\le \limsup_{\varepsilon \to 0} \varepsilon^\alpha [-\log \mathbb{P}(\|X_\rho\|_q < \varepsilon)] \le c_2 \|\rho\|_q. \end{aligned} \qquad (7.9)$$

**Proof.** The left-hand estimate was proven in Li and Linde (2004). To verify the right-hand one, as before, we set $p := q'$ and note that $p > 1$. The operator $u: L_p[0,1] \to L_\alpha[0,1]$ generating $X_\rho$ is given by

$$(uf)(s) = \int_s^1 f(t)\rho(t)\,\mathrm{d}t, \qquad 0 \leq s \leq 1.$$

Because $\rho \in L_{p'}[0,1]$, the corresponding operator $u_\infty$ is well defined from $L_p[0,1]$ into $L_\infty[0,1]$. Upper estimates for the entropy of the operator $u_\infty$ were given in Section 4.6 in Lifshits and Linde (2002). From (4.60) in Lifshits and Linde (2002) (note that $\chi \equiv 1$ and $\eta = \rho$ in the notation of Lifshits and Linde (2002)), we derive

$$\limsup_{n\to\infty} n e_n(u_\infty) \leq c\|\rho\|_{p'}. \tag{7.10}$$

Observe that the right-hand side in (4.58) in Lifshits and Linde (2002) is finite since $r = p'$ and $\chi \equiv 1$. Thus, (4.60) in Lifshits and Linde (2002) applies in our situation and the right-hand estimate in (7.9) follows from Theorem 1.3 using (7.10) and $q = p'$. $\square$

**Remark.** The example of the weighted $\alpha$-stable Lévy motion shows that the application of Theorem 1.3 is limited in some cases. For example, as seen above, it does not apply in the most interesting case $q = \infty$. Here, we have $p = 1$ and then the operator $u_\infty$ is, in general, not compact. Moreover, even if $p > 1$, there remains a gap in the dependence on $\rho$ between the left- and right-hand estimate in Proposition 7.6. Note that $r < q$.

### 7.1.5. The $\alpha$-stable sheet

Finally, we investigate an S$\alpha$S process indexed by $[0,1]^d$ for some $d \geq 1$. If $u$ from $L_p[0,1]^d$ to $L_\alpha[0,1]^d$ is defined by

$$(uf)(s) := \int_{s_1}^1 \cdots \int_{s_d}^1 f(t)\,\mathrm{d}t_d \cdots \mathrm{d}t_1, \qquad s = (s_1,\ldots,s_d),$$

then the generated S$\alpha$S process $Z_\alpha^d$ is usually called a ($d$-dimensional) $\alpha$-stable sheet. Note that for $\alpha = 2$, we obtain the ordinary $d$-dimensional Brownian sheet. An easy transformation gives $e_n(u) = e_n(\bar{u})$, where $\bar{u}$ from $L_p[0,1]^d$ to $L_\alpha[0,1]^d$ is defined by

$$(\bar{u}f)(t) := \int_0^{t_1} \cdots \int_0^{t_d} f(s)\,\mathrm{d}s_d \cdots \mathrm{d}s_1, \qquad t = (t_1,\ldots,t_d).$$

It is known (see Belinsky (1998) and Dunker *et al.* (1999)) that

$$e_n(\bar{u}: L_p[0,1]^d \to L_\infty[0,1]^d) \preceq n^{-1}(\log n)^{d-1/2}$$

whenever $1 < p \leq \infty$. Hence, Theorem 1.3 applies in this case and leads to

$$-\log \mathbb{P}(\|Z_\alpha^d\|_{L_q[0,1]^d} < \varepsilon) \preceq \varepsilon^{-\alpha}\log(1/\varepsilon)^{\alpha(d-1/2)} \qquad \text{for } q \in [1,\infty). \tag{7.11}$$

***Remark.*** Estimate (7.11) is weaker than the best known lower one. Namely, as shown in Li and Linde (2004), we have

$$\varepsilon^{-\alpha}\log(1/\varepsilon)^{\alpha(d-1)} \preceq -\log \mathbb{P}(\|Z_\alpha^d\|_{L_q[0,1]^d} < \varepsilon)$$

for all $q \in [1, \infty]$. Also, the Gaussian case suggests that the exponent of the log-term in (7.11) should be $\alpha(d-1)$, at least for $q < \infty$. Nevertheless, to the authors' knowledge, (7.11) is the first upper estimate for $\phi(Z_\alpha^d, \varepsilon)$.

## 7.2. Sum of maxima-type processes

Let us consider the following random vector. Let $\xi_{n,l}$ be i.i.d. standard S$\alpha$S random variables, $n, l = 1, 2, \ldots$. Let $(\vartheta_n)$ be some decreasing sequence. Consider the random array

$$X := (\vartheta_n \xi_{n,l})_{n=1,2,\ldots; l=1,\ldots,2^n}.$$

We consider $X$ as a random variable in the space $E' := \ell_1(\ell_\infty^{2^n})$, where we use the notation

$$\ell_p(\ell_q^k) := \{z = (z_{n,l}) \mid n = 1, 2, \ldots; l = 1, \ldots, k, \|z\|_{\ell_p(l_q^k)} < \infty\}$$

with the norm given by

$$\|z\|_{\ell_p(\ell_q^k)} := \left[\sum_{n=1}^\infty \left(\sum_{l=1}^k |z_{n,l}|^q\right)^{p/q}\right]^{1/p},$$

with the obvious modification for $p = \infty$ or $q = \infty$.

In our case, we set $E := \ell_\infty(\ell_1^{2^n})$. Now, let $u$ be the diagonal operator from $E = \ell_\infty(\ell_1^{2^n})$ to $\ell_\alpha(\ell_\alpha^{2^n})$ with diagonal $(\vartheta_n)$, that is, $u: (z_{n,l}) \mapsto (\vartheta_n z_{n,l})$. This operator generates the random vector $X \in E'$ with the small deviations

$$\mathbb{P}(\|X\|_{E'} \leq \varepsilon) = \mathbb{P}\left(\sum_{n=1}^\infty \vartheta_n \max_{1 \leq l \leq 2^n} |\xi_{n,l}| \leq \varepsilon\right),$$

which explains the example's name. It was shown in Aurzada and Lifshits (2008) that such probabilities exhibit a critical behavior when the weights are defined by $\vartheta_n = 2^{-n/\gamma} n^{-\beta/\gamma}$ with $\gamma \leq \alpha$. Namely, for $\vartheta_n = 2^{-n/\gamma} n^{-\beta/\gamma}$ with $\gamma < \alpha$, we have

$$-\log \mathbb{P}(\|X\|_{E'} \leq \varepsilon) \approx \varepsilon^{-\gamma} |\log \varepsilon|^{-\beta}. \tag{7.12}$$

For the 'critical' case $\gamma = \alpha$, that is, $\vartheta_n = 2^{-n/\alpha} n^{-\beta/\alpha}$, however, we have

$$-\log \mathbb{P}(\|X\|_{E'} \leq \varepsilon) \begin{cases} = \infty, & \beta \leq \max(1, \alpha), \\ \approx \varepsilon^{-1/(\beta/\alpha - 1)}, & \max(1, \alpha) < \beta < 1 + \alpha, \\ \approx \varepsilon^{-\alpha} |\log \varepsilon|^{1+\alpha}, & \beta = 1 + \alpha, \\ \approx \varepsilon^{-\alpha} |\log \varepsilon|^{-\beta + 1 + \alpha}, & \beta > 1 + \alpha. \end{cases} \tag{7.13}$$

The entropy numbers of the generating operator $u$, as calculated by T. Kühn (private communication), are

$$e_n(u) \approx n^{1/\alpha - 1/\gamma - 1}(\log n)^{-\beta/\gamma},$$

which agrees with (7.12) from the point of view of Theorem 1.3. However, the connection between the entropy and small deviations completely breaks down in the critical case (7.13).

This example shows that the connection between small deviations and the entropy of the related operator can deviate quite drastically from what would be expected from the valid implications in Proposition 1.2. It is an open problem to calculate a corresponding operator $u_\infty$ in this case and to see to which bounds this leads to by the use of Theorem 1.3.

### 7.3. Small deviation of S$\alpha$S-vectors with $0 < \alpha < 1$

Finally, let us indicate a relation of our result to a general lower bound due to Ryznar (1986). Namely, note that, trivially,

$$e_n(u_\infty) = e_n(u : E \to L_\infty(S)) \leq \|u : E \to L_\infty(S)\| < \infty.$$

Therefore, we can set $\tau := \alpha/(1-\alpha)$ (which is positive for $0 < \alpha < 1$) and $L = 1$ in Theorem 1.3 and obtain that, for any S$\alpha$S random vector with $0 < \alpha < 1$,

$$\phi(X, \varepsilon) \preceq \varepsilon^{-\alpha/(1-\alpha)}.$$

This result was shown in Ryznar (1986) for strictly stable (not necessarily symmetric) vectors with $0 < \alpha < 1$, using completely different methods.

## Acknowledgements

We would like to thank J. Christof (Jena) for reading the manuscript carefully and pointing out several misprints. The research of the first author was supported by the DFG Research Center MATHEON "Mathematics for key technologies" in Berlin. The joint work of the second and the third authors was supported by the DFG-RFBR Grant 09-01-91331.